\documentclass[12pt]{article}

\usepackage[T1]{fontenc}
\usepackage{lmodern} 
\usepackage{amsmath,amsfonts,amssymb,amsthm}
\usepackage{graphicx}
\usepackage[usenames,dvipsnames,table]{xcolor}
\usepackage{xspace}
\usepackage[normalem]{ulem}
\usepackage{subcaption}
\usepackage{tikz}
\usepackage[inkscapeformat=pdf]{svg}

\usepackage[letterpaper,top=2.8cm,bottom=2cm,left=2.9cm,right=2.9cm,marginparwidth=2cm]{geometry}
\usepackage[colorlinks=true, allcolors=blue, breaklinks=true]{hyperref}

\usepackage{booktabs}
\usepackage{siunitx}

\usepackage{authblk}

\usepackage{natbib}

\definecolor{darkviolet}{rgb}{0.58, 0.0, 0.83}

\newcommand{\mlrule}{\texttt{ML-Q-RF}\xspace}
\newcommand{\optimal}{\texttt{ORule\textsuperscript{S}}\xspace}
\newcommand{\bestrule}{\texttt{BRule\textsuperscript{D}}\xspace}
\newcommand{\bestvar}{\texttt{BVar\textsuperscript{D}}\xspace}
\newcommand{\bestruleQRF}{\texttt{BRule\textsuperscript{D,Q-RF}}\xspace}
\newcommand{\bestvarQRF}{\texttt{BVar\textsuperscript{D,Q-RF}}\xspace}
\newcommand{\bestruleOPT}{\texttt{BRule\textsuperscript{D,opt}}\xspace}
\newcommand{\bestvarOPT}{\texttt{BVar\textsuperscript{D,opt}}\xspace}
\newcommand{\bestvarOPTtie}[1]{\texttt{BVar\rlap {\textsubscript{\tiny{#1}}}\textsuperscript{D,opt}}\xspace}
\newcommand{\bestvarQRFtie}[1]{\texttt{BVar\rlap {\textsubscript{\tiny{#1}}}\textsuperscript{D,Q-RF}}\xspace}

\usepackage{todonotes}

\title{Learning in Spatial Branching: Limitations of Strong Branching Imitation}
\author[1]{Brais González-Rodríguez}
\author[1,2]{Ignacio Gómez-Casares}
\author[3]{Bissan Ghaddar}
\author[1,2]{Julio González-Díaz}
\author[1,2]{Beatriz Pateiro-L\'opez}
\affil[1]{Department of Statistics, Mathematical Analysis and Optimization and MODESTYA Research Group, University of Santiago de Compostela, 15782 Santiago de Compostela, Spain}
\affil[2]{CITMAga (Galician Center for Mathematical Research and Technology), 15782 Santiago de Compostela, Spain}
\affil[3]{Ivey Business School, Western University, London, Ontario, Canada}

\date{}

\begin{document}

\maketitle

\begin{abstract}
Over the last few years, there has been a surge in the use of learning techniques to improve the performance of optimization algorithms. In particular, the learning of branching rules in mixed integer linear programming has received a lot of attention, with most methodologies based on strong branching imitation. Recently, some advances have been made as well in the context of nonlinear programming, with some methodologies focusing on learning to select the best branching rule among a predefined set of rules leading to promising results. In this paper we explore, in the nonlinear setting, the limits on the improvements that might be achieved by the above two approaches: learning to select the best variable (strong branching) and learning to select the best rule (rule selection). 
\end{abstract}



\section{Introduction}\label{sec:intro}
{In the field of optimization and, more specifically, in the context of branch-and-bound algorithms, branching decisions significantly affect the efficiency and speed of finding the optimal solution. Strong branching is a specific type of branching technique that involves evaluating all possible branching variables by simulating the branching process for each potential variable and measuring the impact on the optimality gap. Although strong branching is known to be highly effective in reducing the size of the search tree and accelerating the convergence to an optimal solution, it is computationally prohibitive because it requires solving, at each and every node of the branch-and-bound tree, two optimization problems for each candidate variable. Learning to branch leverages machine learning techniques to imitate or enhance traditional branching strategies like strong branching. By using historical data on past branching decisions and their outcomes, machine learning models are trained to approximate the decision-making quality of strong branching without incurring high computational costs, thus balancing efficiency and computational speed.} As can be seen in the recent surveys \cite{Lodi2017}, \cite{Bengio2021}, \cite{fan2024}, and \cite{scavuzzo2024}, most of the research so far has centered on mixed integer linear programming due to its wide range of applications. This paper contributes to the same strand of literature, but in the context of nonlinear problems.

When it comes to the comparison of learning methodologies, one can establish different taxonomies. An important distinction is between offline learning and online learning. In the context of branching, offline learning requires that the training is separated from the solution of a new optimization problem: the (possibly) computationally demanding training phase takes place upfront, and there is no computational overhead at the moment of execution of the optimization algorithm. On the other hand, online learning requires gathering data node by node as the branch-and-bound tree grows, and the learning is dynamically updated as new data becomes available: the resulting branching rules are more adaptive, but there is a trade-off between the potential improvements in performance and the computational cost of learning, which is now part of the cost of execution. As far as features are concerned, one can distinguish between static and dynamic features. Static features only depend on the specific characteristics of each optimization problem, such as the number of variables, the number of constraints, some descriptive statistics on the problem coefficients, or the sparsity of the problem. Dynamic features, on the other hand, can be much richer, capturing information about the evolution of the branch-and-bound algorithm such as the progress in the bounds of the problem, bounds on the variables, and even properties of the branch-and-bound tree itself. 


Returning to the specific learning task of branching variable selection, a common approach involves using a known expert system as a model and training on past data to imitate it. As discussed above, given that strong branching is (experimentally) known to produce significantly smaller trees than any other branching rule, it is not surprising that, in this context, most literature has focused precisely on learning to imitate strong branching, as in \cite{Khalil2016}, \cite{Alvarez2017}, \cite{Balcan2018}, \cite{Gasse2019}, and \cite{gupta2020}. {For instance, \cite{Alvarez2017} present an offline supervised learning method to predict strong branching scores in branch-and-bound algorithms for solving mixed-integer linear programming problems. \cite{Khalil2016} explore an online framework that ranks variables to replicate strong branching. On the other hand, \cite{Gasse2019} use offline learning to encode the branching policies into a graph convolutional neural network and have fast approximations for strong branching. \cite{gupta2020} extend this work to a hybrid offline learning model that combines the capabilities of graph neural networks with the less computationally demanding multi-layer perceptrons for more efficient branching. }

Despite the widely spread use of strong branching as an expert to imitate, recent research has highlighted some limitations associated with this approach. As discussed in \cite{scavuzzo2024}, \emph{``The strong branching heuristic has been used by many as an expert from which we can learn effective decision-making. The claim that strong branching is a desirable strategy to follow has recently been challenged, with some notable examples of instances where strong branching scores provide no useful information''.} \cite{gamrath2018} argue that, when measuring the impact of strong branching on metrics such as the number of nodes of the branch-and-bound tree, careful consideration is needed to obtain fair comparisons. In particular, since strong branching can identify infeasibilities (and does not perform branching at nodes where it identifies such infeasibilities) and provide faster dual-bound improvements, standard implementations of strong branching effectively reduce the size of the resulting tree in a way that no rule aimed at imitating strong branching can ever achieve. Second, as also highlighted in \cite{scavuzzo2024}, \emph{``strong branching relies on the LP relaxation for scoring variables, which can provide little information in cases when the optimal LP objective value does not change with branching. In such cases, strong branching cannot be considered a reliable expert.''} These concerns have motivated research on expert-free approaches, in which reinforcement learning is used to learn a branching rule from scratch as in \cite{etheve2020}, \cite{Parsonson2022}, and \cite{scavuzzo2022}. 


Although the use of learning in nonlinear programming (NLP) has also been on the rise in recent years (see, for instance, \cite{Baltean2019}, \cite{bonami2022}, and \cite{kannan2023}), not much has been done on learning to branch. An exception is \cite{ghaddar2023}, where the authors develop an offline learning framework for spatial branching based on ideas of algorithm selection \citep{Gomes2001}, similar to those used in \cite{Liberto2016} for MILP problems. In this paper we extend the work in \cite{ghaddar2023}, in which learning takes place at the instance level with static features only, and study to what extent learning at the node level with dynamic features can outperform the baseline static approach. On the one hand, for a given set of optimization problems, the dynamic setting provides many more observations for the training. Each instance results in as many observations as nodes of the branch-and-bound tree, instead of having just one observation associated with each instance. On the other hand, the information obtained on a node-by-node basis is much noisier, which can definitely be challenging for the learning.

The contribution of this paper is not to assess the performance of a specific learning technique, but to determine the potential of a learning framework itself. More specifically, we evaluate the potential of three learning frameworks, by analyzing the performance of three different experts. The first one, \optimal, is an expert capable of choosing, for each instance, the best performing branching rule from a predefined pool; it represents the best that can be achieved with the static learning in \cite{ghaddar2023}. The second one, \bestrule, is an expert capable of choosing, for each node of the branch-and-bound tree, the best performing branching rule (at that node) from a predefined pool; it represents the best that can be achieved by extending the static learning setting to accommodate dynamic features and learn at the node level. The third one, \bestvar, is an expert capable of choosing, for each node of the branch-and-bound tree, the best variable (at that node), in the sense of perfectly imitating strong branching choices. Importantly, note that experts \bestrule and \bestvar can be seen as targets both for offline and online learning and, indeed, as we discussed above, \bestvar has been the most common choice by past literature under both settings. Therefore, our findings apply to both offline and online learning frameworks.


We frame our contribution in the context of continuous polynomial optimization problems, where we study branching variable selection for a spatial branch-and-bound algorithm that builds upon the RLT technique \citep{Sherali1992}. Our numerical results show that \optimal is the best performing expert: both \bestrule and \bestvar fall behind for the various performance metrics under consideration. This finding might seem counter intuitive at first, since the last two experts should benefit from the richness of the dynamic setting and also from the possibilities of using online learning or hybrid offline-online approaches. Yet, what seems to be limiting these experts is the fact that they make, node by node, the best myopic branching decision based on the achieved lower bound improvement at the current node. Thus, there is no guarantee that this locally optimal choice will result in a good enough overall performance and, indeed, we observe that \optimal, the expert of the (offline) static framework is clearly superior.

It is worth noting that the use of an RLT-based algorithm for polynomial optimization problems is merely a means of assessing the potential of the different experts and learning frameworks. Thus, although polynomial optimization is an important class of problems on its own right, we believe that the implications and insights of our findings can be extrapolated to general NLP and MINLP problems. In particular, the relatively weak performance obtained for \bestvar, might guide future efforts in the design of branching rules for spatial branching in directions whose goal is not necessarily to approximate strong branching decisions.

The rest of the paper is structured as follows. In Section~\ref{sec:framework} we present the overall setting and in Section~\ref{sec:testing} we discuss the testing environment. Section~\ref{sec:results} contains the numerical analysis. We conclude in Section~\ref{sec:conclusions}.

\section{Framework for the analysis}\label{sec:framework}
In this section, we describe the elements required to develop the computational study to assess the potential of learning schemes based on the imitation of three different experts: \optimal, \bestrule, and \bestvar.

\subsection{Optimization problem and branch-and-bound algorithm}
We contextualize our analysis within the class of polynomial optimization problems and global optimization schemes based on the use of spatial branching based on reformulation-linearization technique. Formally, we consider (continuous) polynomial optimization problems of the following form:
\begin{equation}
\begin{split}
\text{minimize} & \quad \phi_0(\mathbf{x})\\
\text{subject to}  & \quad \phi_r(\mathbf{x})\geq \beta_r, \quad r=1,2,\ldots, R_1 \\
& \quad \phi_r(\mathbf{x})=\beta_r, \quad r=R_1+1,\ldots,R\\
& \quad \mathbf{x}\in\Omega \subset \mathbb{R}^{|N|}\text{,}
\end{split}
\label{eq:PO}
\tag{PO}
\end{equation}
where $N$ denotes the set of variables, each $\phi_r(\mathbf{x})$ is a polynomial of degree $\delta_r \in \mathbb{N}$, and the region $\Omega = \lbrace \mathbf{x} \in \mathbb{R}^{|N|}: 0 \leq l_j \leq x_j \leq u_j < \infty, \, \forall j \in N \rbrace \subset \mathbb{R}^{|N|}$ is a hyperrectangle containing the feasible region. Then, $\delta=\max_{r \in \{0,\ldots,R\}} \delta_r$ is the degree of the problem and $(N, \delta)$ represents all possible monomials of degree $\delta$.

The \textit{reformulation-linearization technique}, RLT, was introduced in \cite{Sherali1992} to solve polynomial optimization problems like~\eqref{eq:PO}. This is achieved by defining a linear relaxation of~\eqref{eq:PO}, which is then embedded in a branch-and-bound scheme. Each monomial of degree greater than one is replaced with a corresponding RLT variable. For example, the RLT variable $X_{1224}$ would replace  monomial $x_1x_2x_2x_4$. More generally, RLT variables are defined as
\begin{equation}
  X_J = \prod_{j \in J}x_j, 
  \label{eq:RLTidentity}
\end{equation}
where $J$ is a multiset containing the information about the multiplicity of each variable in the underlying monomial. These relaxations are then solved at each node of the branch-and-bound tree. In order to get tighter relaxations and ensure convergence of the resulting algorithm, additional constraints, called bound-factor constraints, are also added {(and linearized)}. They are given, for each pair of multisets $J_1$ and $J_2$ such that $J_1\cup J_2 \in (N,\delta)$, by $F_{\delta}(J_1,J_2)=\prod_{j\in J_1}{(x_j-l_j)}\prod_{j\in J_2}{(u_j-x_j)}\geq 0$.

\subsection{Porfolio of branching rules}\label{sec:rules}
We next describe how the branching decisions, which are the focus of our study, is carried out in the RLT-based algorithm. Upon solving the corresponding relaxation at each node, a score is assigned to each variable based on the violations of RLT identities~\eqref{eq:RLTidentity} in which it is involved. Then, the variable with the highest score is chosen for branching.

State-of-the-art solvers for non-linear optimization rely on two main themes for branching variable selection. The first one involves adaptations of MILP methods, such as strong branching and reliability branching \citep{Achterberg2005}, as explained in \cite{Belotti2009}. The second one consists of identifying variables that contribute to the violations of non-linear terms, using variations of the violation transfer method, which was introduced in \cite{Tawarmalani2004} and further developed in \cite{Belotti2009}. 
In \cite{Gonzalez-Rodriguez:2023} and \cite{ghaddar2023} these approaches were studied in conjunction with the violations of the RLT-defining identities, which have traditionally been the main metric for variable selection in the RLT literature. We briefly describe below the best performing approaches in \cite{ghaddar2023}, which we take as the starting point in our analysis. All the approaches are of the form
\begin{equation}
\label{eq:sumweightscriteria}
  \theta_j = \sum_{J\subset (N,\delta)\colon |J|<\delta} w(j,J) \vert \bar{X}_{J \cup \lbrace j \rbrace} - \bar{x}_j \bar{X}_{J} \vert,
\end{equation}
where $w(j,J)$ are the weights associated with the violations and $J \cup \{j\}$ considers only the monomials present in~\eqref{eq:PO}. Building upon this formula we consider the following six branching rules:

\begin{description}
	\item[\texttt{Dual}.] Eq.~\eqref{eq:sumweightscriteria}, with $w(j,J)$ defined as the sum of the absolute values of the shadow prices of the constraints containing $J\cup \{j\}$, at the optimal solution of the current node's relaxation.
	\item[\texttt{Range}.] Eq.~\eqref{eq:sumweightscriteria}, with $w(j,J)=\frac{\min{\lbrace\bar{u}_j - \bar{x}_j, \,\bar{x}_j - \bar{l}_j\rbrace}}{u_j-l_j}$, where $u_j$ and $l_j$ are the upper and lower bound of $x_j$ at the root node and $\bar{u}_j$ and $\bar{l}_j$ are the bounds at the current node. 
	\item[\texttt{Eigen}.]Eq.~\eqref{eq:sumweightscriteria}, with $w(j,J)$ defined as the eigencentrality of $x_j$'s node in the {constraints-monomials intersection graph} associated to \eqref{eq:PO}.
	\item[\texttt{Dual Rel}, \texttt{Range Rel}, \texttt{Eigen Rel}.] These three additional branching rules are defined by embedding reliability branching ideas \citep{Achterberg2005,Belotti2009} into the weight computation of the three rules above.
\end{description}
Once we have defined the portfolio of branching rules, we move to the learning frameworks and associated experts.


\subsection{Learning frameworks and experts}\label{subsec:experts}
In this section, we formally present the three learning frameworks and the three associated experts whose potential is the main object of analysis in this paper. Each of them represents, in a sense that we clarify below, the best that might ever be achieved within a certain learning framework. Therefore, as already mentioned, the contribution of this paper is not to assess the performance of a specific learning technique, but to determine the potential of a learning framework itself. More precisely, we want to understand how much a branch-and-bound algorithm for nonlinear optimization might possibly benefit by enhancing branching variable selection with machine learning. This question is studied for three different learning frameworks, which we describe below.

\subsubsection{Instance-by-instance rule learning with static features}\label{sec:staticrule}

Consider an offline learning setting in which there is a pool of branching rules available and the goal is to learn to choose, instance by instance, the best performing branching rule according to a given metric. Since the learning is performed at the instance level, it can only use features that depend on the specific characteristics of each optimization problem, \emph{i.e.}, static features. Hereafter, we refer to this learning framework as \emph{static rule learning}.

The absolute best that can be achieved with any learning scheme in the static learning setting is to perfectly predict the best performing branching rule on all instances. \optimal represents an expert capable of making these perfect predictions. Note that this expert is optimal in the following sense: no learning scheme designed to assign a branching rule from the pool to each instance can ever outperform \optimal. 

A limitation of the static learning framework is that the learning is not adaptive, since one commits to applying the same branching rule at all the nodes of the branch-and-bound tree.

\subsubsection{Node-by-node rule learning with dynamic features}\label{sec:dynamicrule}
A natural extension of the above setting is to carry out the learning scheme on a node-by-node basis. This approach allows the incorporation of dynamic features into the learning process using, for instance, information about past branching decisions, the evolution of the bounds, or the structure of the branch-and-bound tree. In order to measure the performance at a given node, we take as \emph{key performance indicator}, KPI, the improvement with respect to the lower bound of the parent node. Therefore, the goal is to learn to choose, node by node, the branching rule from the pool that will deliver the best KPI.  We refer to this framework as \emph{dynamic rule learning}.

We use \bestrule to refer to an expert capable of perfectly predicting the branching rule with the best KPI at each node. Note that, for this expert, we have chosen to refer to it as ``Best''  instead of ``Optimal'', because the chosen KPI is a myopic one, since it is only a measure of the improvement at the current node, \emph{i.e.}, it does not capture the performance further down the tree. Expert \bestrule can be seen as the target to imitate in any setting with node-by-node learning, regardless of whether the learning takes place offline or online.

Despite the advantages of the additional richness of the dynamic setting, node-by-node information can also be very noisy, which poses new challenges. Further, it is uncertain whether the optimal node-by-node decisions of \bestrule will result in a good overall performance, since they might be limited by the myopic nature of the KPI. This potential limitation of \bestrule is also shared by the next expert, strong branching imitation. Nonetheless, the latter is still the most common approach in the learning-to-branch literature (as discussed in Section \ref{sec:intro}).

\subsubsection{Node-by-node variable learning with dynamic features}\label{sec:dynamicvariable}
The two experts defined above, \optimal and \bestrule, focus on learning to select the best rule from a given pool or branching rules. Yet, most of the research in the field focuses instead on learning the best branching variable and, more precisely, on strong branching imitation, which is the approach behind our third expert. Consider again the learning setting with dynamic features and the KPI used to introduce \bestrule. The goal now is to choose, node by node, the branching variable that provides the largest improvement with respect to the lower bound of the parent node. We refer to this framework as \emph{dynamic variable learning}.

We use \bestvar to refer to an expert capable of perfectly predicting the branching variable with the best KPI at each node. Note that, as defined, \bestvar coincides with (perfect) strong branching imitation. The above discussion for \bestrule on the use of ``Best'' instead of ``Optimal'' and on its potential limitations because of the myopia of the KPI also apply to \bestvar. Expert \bestvar can be seen as the target to imitate in any setting with node-by-node learning, regardless of whether the learning takes place offline or online.

\section{Testing environment}\label{sec:testing}

For the numerical results presented in this paper, we use the polynomial optimization solver RAPOSa 3.0.1 \citep{Gonzalez-Rodriguez:2023}, whose core is an RLT-based algorithm.\footnote{RAPOSa was run using J-sets \citep{Sherali2013}, while other enhancements such as warm start, bound-tightening techniques \citep{Belotti:2012, gomezcasares2024} and SDP-cuts and conic relaxations \citep{Sherali2012cuts, raposaconic} were disabled.} For the computational experiments, we use a superset of instances that borrows from three well-known test sets. The first one is taken from \cite{Dalkiran2016} and consists of 180 instances of randomly generated polynomial optimization problems of different degrees, number of variables, and density. The second test set comes from MINLPLib \citep{minlplib}, a reference library of mixed-integer nonlinear programming instances. We have selected from MINLPLib those instances that are polynomial optimization problems with box-constrained and continuous variables, resulting in a total of 166 instances. The third test set comes from another well-established benchmark, QPLIB \citep{qplib}, a library of quadratic programming instances, for which we made a selection analogous to the one made for MINLPLib, resulting in a total of 63 instances. 
All the executions reported in this paper have been run on the supercomputer Finisterrae~III, at Galicia Supercomputing Centre (CESGA). Specifically, we used computational nodes powered with two thirty-two-core Intel Xeon Ice Lake 8352Y CPUs with 256GB of RAM and 1TB SSD. All the executions were been run with a time limit of 10 minutes and the stopping criterion was set at threshold $0.001$ for either the relative gap or the absolute gap.

To compare the performance of the different approaches under consideration, we use the following metrics:
\begin{itemize}
 \item \textbf{Solved}. Number of instances solved to certified optimality (relative or absolute gap below $0.001$).
 \item \textbf{Gap}. Geometric mean of the optimality gap obtained by each approach. We exclude instances for which i)~at least one approach did not return an optimality gap and ii)~all the approaches solved it within the time limit.
 \item \textbf{Time}. Geometric mean of the running time of each approach. We exclude instances for which i)~every approach solved it in less than 5 seconds and ii)~no approach solved it within the time limit.
 \item \textbf{Pace}. Geometric mean of pace$^{LB}$, as introduced in \cite{ghaddar2023}. It represents the number of seconds needed to improve the lower bound of the algorithm by one unit (the pace at which the lower bound improved). We exclude instances solved by every approach in less than 5 seconds. The main motivation behind this performance measure is that it allows to compare the performance on all instances together, whereas Time and Gap fail to do so. Time is not informative when comparing performance between instances not solved by any approach (all of them reach the time limit). These instances can be compared using Gap which, in turn, is not informative in instances solved by all approaches (all of them close the gap), and where Time could be used. Pace, on the other hand, is informative regardless of the number of approaches that might have solved each of the different instances. 
 \item \textbf{Nodes}. Geometric mean of the number of nodes explored in the branch-and-bound tree. We only consider the instances solved by all the approaches within the time limit.
\end{itemize}

The metric-specific exclusions described above facilitate the interpretation of the results associated with the different metrics. These exclusions imply that, in each table of results, the precise numbers describing the performance of each approach depend on the whole set of approaches being compared, since all of them are used to determine the exclusions for each performance metric. The number of instances used for the computation of each performance metric is shown in parentheses. 


\section{Assessing learning frameworks}\label{sec:results}
In this section, we present the computational analysis of the three learning frameworks discussed in Section~\ref{subsec:experts}: static rule learning, dynamic rule learning, and dynamic variable learning. In order to analyze the potential of these three learning frameworks, we study the performance of the associated experts: \optimal, \bestrule, and \bestvar, respectively.

\subsection{Static rule learning}

We take as our starting point the learning framework developed in \cite{ghaddar2023} in the context of what we have called static rule learning. The KPI used in their analysis is a normalized version of pace$^{LB}$, which allows training on a set of instances of varying degrees of difficulty, for which neither time nor gap would be suitable KPIs. The improvements reported there are substantial, with the machine learning version delivering improvements of up to $25-35\%$ with respect to the best original branching rule, whereas $\optimal$ would lead to an additional $15-20\%$ improvement. Table~\ref{tab:static-orig} reports the results of the application of their learning methodology to the pool of rules comprised of the six branching rules defined in Section~\ref{sec:rules}. In particular, the learning is conducted jointly on all sets of instances using quantile regression forests \citep{Meinshausen2006}, and the rule resulting from the learning process is referred to as \mlrule.\footnote{The statistical analysis was developed in programming language~\texttt{R} \citep{Rlang}, using library \texttt{ranger}~\citep{rangerR}.}

\sisetup{
	round-mode = places,
	round-precision = 2,
	detect-weight=true,
	table-number-alignment = right,
	table-alignment-mode = none,
	table-text-alignment = right
	}
 
\begin{table}[!htpb]
    \caption{Pool of branching rules and resulting machine learning rule.}
    \label{tab:static-orig}
    \centering
    {\addtolength{\tabcolsep}{-2pt} 
    \begin{tabular}{l| r S[round-precision=3] S S S}
    \toprule
        {} & {Solved {\scriptsize(409)}} & {Gap {\scriptsize(159)}} & {Time {\scriptsize(120)}} & {Pace {\scriptsize(265)}} & {Nodes {\scriptsize(235)}} \\ \hline
        {\texttt{Dual}} & 254 & 0.113938 & 46.022532 & 9.338496 & 98.991928\\
        {\texttt{Range}} & 249 & 0.113895 & 51.196735 & 9.764979 & 94.739696\\
        {\texttt{Eigen}} & 243 & 0.135928 & 93.852929 & 11.210842 & 135.711555\\
        {\texttt{Dual Rel}} & 249 & 0.130042 & 45.245676 & 9.455919 & 89.520644 \\
        {\texttt{Range Rel}} & 252 & 0.115264 & 46.253278 & 6.731057 & 88.889385\\
        {\texttt{Eigen Rel}} & 242 & 0.144880 & 78.118226 & 10.225294 & 118.092436\\
        {\texttt{ML-Q-RF}} & 257 & 0.096719 & 39.702132 & 4.563422 & 89.099049 \\
    \bottomrule
    \end{tabular}}
\end{table}


The results in Table~\ref{tab:static-orig} show that \mlrule improves upon the best performing rule according to Pace, \texttt{Range Rel}, by $32\%$ (from $6.73$ to $4.56$) and delivers also significant improvements in Solved, Time, and Gap.\footnote{Our results do not reproduce exactly those reported in \cite{ghaddar2023}, since there are some changes in the computational setup, such as the version of RAPOSa used in the analysis.}  In Table~\ref{tab:static-optimal} we assess the performance of our first expert, \optimal, by comparing it with \texttt{Range Rel} and \mlrule.\footnote{Recall that the results in Table~\ref{tab:static-optimal} for \texttt{Range Rel} and \mlrule do not need to coincide with those in Table~\ref{tab:static-orig}. This is because of the exclusions applied to compute the values for Gap, Time, Pace, and Nodes described in Section~\ref{sec:testing}, which leads to slightly different numbers of instances considered for each metric.} In this case, the difference in Pace between \texttt{Range Rel} and \mlrule is 33\%, whereas \optimal delivers an additional $14\%$, almost reducing Pace to half (from $9.19$ to $4.87$) and producing substantial improvements on the five performance metrics under study.

\begin{table}[!htpb]
    \caption{Assessing the potential of \optimal.}
    \label{tab:static-optimal}
    \centering
    {\addtolength{\tabcolsep}{-2pt} 
    \begin{tabular}{l| r S[round-precision=3] S S S}
    \toprule
        {} & {Solved {\scriptsize(409)}} & {Gap {\scriptsize(145)}} & {Time {\scriptsize(108)}} & {Pace {\scriptsize(253)}} & {Nodes {\scriptsize(250)}} \\ \hline
        {\texttt{Range rel}} & 252 & 0.182285 & 69.775062 & 9.188080 & 115.472324\\
        {\texttt{ML-Q-RF}} & 257 & 0.150388 & 59.565470 & 6.145497 & 108.254329\\
        \optimal & 259 & 0.134097 & 45.575315 & 4.870658 & 92.763685 \\
    \bottomrule
    \end{tabular}}
\end{table}


Given the potential of both \mlrule and \optimal illustrated in Tables~\ref{tab:static-orig} and \ref{tab:static-optimal}, one may want to explore to what extent the results could be further improved by enriching the relatively simple static rule learning framework above. A first approach would be to add more diverse and competitive branching rules to the pool since, by construction, the more rules in the pool, the better \optimal will perform. This idea was already explored in \cite{ghaddar2023}, where the authors found that \mlrule also benefits from the additional richness. A second approach, which we consider here, consists of adapting the algorithm selection approach to the dynamic rule learning framework, \emph{i.e.}, learning to select the best rule on a node-by-node basis using dynamic features.

\subsection{Dynamic rule learning}
As discussed, the promising results of the previous section have been obtained in a setting in which features are restricted to be static and where the size of the training set is limited by the available number of instances. Thus, it is reasonable to question how much further one might go in a setting where learning takes place at the node level, where new features would be available to capture the impact of past decisions and the evolution of the branch-and-bound algorithm {and where, moreover, one might also rely on online learning}. This is what we referred to as dynamic rule learning in Section~\ref{sec:dynamicrule}, and \bestrule is an expert capable of perfectly predicting the branching rule with the best KPI at each node.

It is important to note that there will be nodes along the branch-and-bound tree for which different rules select different branching variables and, still, have the same KPI. Probably the most common (and relevant) situation in which this happens is when all rules lead to a KPI of 0, \emph{i.e.}, no improvement with respect to the lower bound of the parent node. Therefore, in order to perfectly define \bestrule, we have to specify how to proceed when such ties arise. Since the main goal of this section is to assess the extra potential of the dynamic rule framework with respect to the static one, we explore two possibilities, \bestruleQRF and \bestruleOPT, which in case of a tie follow, respectively, the choice that \mlrule and \optimal would make for the instance under consideration. In particular, \bestruleOPT allows the expert of the dynamic setting to build upon the expert of the static learning framework {in situations where the chosen KPI produces ties}. 

\begin{table}[!htpb]
    \caption{Assessing the potential of \bestrule.}
    \label{tab:dynamic-rule}
    \centering
    {\addtolength{\tabcolsep}{-2pt} 
    \begin{tabular}{l| r S[round-precision=3] S S S}
    \toprule
        {} & {Solved {\scriptsize(409)}} & {Gap {\scriptsize(142)}} & {Time {\scriptsize(109)}} & {Pace {\scriptsize(252)}} & {Nodes {\scriptsize(256)}} \\ \hline
        {\texttt{ML-Q-RF}} & 257 & 0.191486 & 63.340109 & 6.478339 & 125.215455\\
        \optimal & 259 & 0.169800 & 48.356430 & 5.119357 & 105.144451\\
        \bestruleQRF & 261 & 0.183248 & 53.784458 & 6.376590 & 104.506109\\
        \bestruleOPT & 260 & 0.193972 & 53.632602 & 6.132581 & 105.192220\\
    \bottomrule
    \end{tabular}}
\end{table}


Table~\ref{tab:dynamic-rule} presents the performance of \bestruleQRF and \bestruleOPT side by side with the performance of \mlrule and \optimal. Surprisingly, we see that \optimal is significantly superior to both \bestruleQRF and \bestruleOPT in terms of Pace, Gap, and Time, while falling narrowly behind in the number of solved instances. Thus, it seems that the myopic nature of the chosen KPI represents a ceiling on what can be achieved in the dynamic rule framework, whose potential is below what can be achieved by \optimal. Interestingly, even \mlrule is highly competitive with both variants of \bestrule, and specially so in terms of Pace, the KPI on which it was trained. These findings suggest that it may be more promising to enhance the relatively simple approaches in static rule learning, rather than pursuing what are likely to be much more complex techniques in the context of dynamic rule learning. 
  
Figure~\ref{fig:ppbestvar} presents three performance profiles \citep{Dolan2002}, which provide additional information to compare \mlrule, \optimal, and \bestruleOPT, according to running time, optimality gap, and pace. The three of them are consistent with the information in Table~\ref{tab:dynamic-rule}, but the superiority of \optimal is even more noticeable.

\begin{figure}[!htpb]
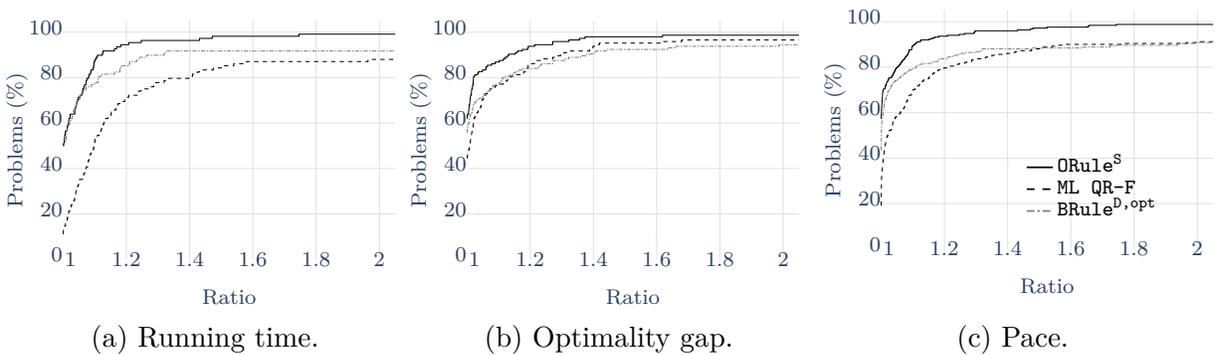

\centering
\begin{subfigure}{0.32\textwidth}
\includesvg[pretex=\scriptsize,width=\textwidth]{images/time1.svg}
    \caption{Running time.}
    \label{subfig:tim1}
\end{subfigure}
\hfill
\begin{subfigure}{0.32\textwidth}
    \includesvg[pretex=\scriptsize,width=\textwidth]{images/gap1.svg}
    \caption{Optimality gap.}
    \label{subfig:gap1}
\end{subfigure}
\hfill
\begin{subfigure}{0.32\textwidth}
    \begin{tikzpicture}
     \node (fig1) at (0,0){\includesvg[pretex=\scriptsize,width=\textwidth]{images/pace1.svg}};
     \node (fig2) at (1,-0.4)
       {\includesvg[pretex=\scriptsize,width=0.4\textwidth]{images/legend1_JGD.svg}};  
    \end{tikzpicture}
    \caption{Pace.}
    \label{fig:pace1}
\end{subfigure}       
\caption{Performance profiles.}
\label{subfig:ppbestrule}
\end{figure}

\subsection{Dynamic variable learning}
Given the negative results obtained for \bestrule in the previous section, it is important to get a better understanding on why the dynamic rule learning framework has a more limited potential than the much simpler static one. As discussed above, a natural explanation is that the underlying KPI is completely myopic, so even an expert capable of following the rule with the best KPI at each node may not exhibit a good overall performance. Alternatively, one may argue that the limitation comes from the algorithm selection approach, under which the only candidate variables for branching at a given node are those chosen by at least one of the branching rules from the pool. In the latter case, we should observe that an expert capable of choosing, at each node, the branching variable with the best KPI, should outperform \bestrule. This is exactly what we study in this section. Evaluating the performance of such an expert, which we have called \bestvar, is equivalent to evaluating the performance of perfectly imitating strong branching, which, as discussed in the introduction, has been the goal of most of the literature on learning to branch in MILP problems.

Similar to the above section with \bestrule, we also need to specify what variables are selected for branching when ties in KPI occur. Again, we consider two variants of \bestvar: \bestvarQRF and \bestvarOPT, which in case of a tie follow, respectively, the choice that \mlrule and \optimal would make for the instance under consideration. There is another aspect regarding ties that will be relevant in the ensuing analysis, namely, whether or not we should include some tolerance when considering ties. To this end, \bestvarQRFtie{0} and \bestvarOPTtie{0} denote the versions of \bestruleQRF and \bestruleOPT in which only a perfect coincidence of the KPI (to machine accuracy) is considered a tie. On the other hand, we later consider as well variants such as \bestvarOPTtie{0.005} or \bestvarOPTtie{0.01}, to follow \optimal also in situations where the best KPI delivers an improvement no superior to $0.005$ or $0.01$, respectively.

Given that strong branching is computationally very intensive, and that this burden rapidly increases with the number of variables, for the results reported in this section we have restricted our test set to instances with 20 variables or less. This reduces the size of the test set from 409 instances to 240. It is worth noting that, for the numeric results below, we assume that \bestvar and its variants can perfectly imitate strong branching with no computational overhead, \emph{i.e.}, they are just assumed to know, for each and every node of the branch-and-bound tree, the KPI associated to each and every variable.

\begin{table}[!htpb]
    \caption{Assessing the potential of \bestvar.}
    \label{tab:dynamic-sb00}
    \centering
    {\addtolength{\tabcolsep}{-2pt} 
    \begin{tabular}{l| r S[round-precision=3] S S S}
    \toprule
        {} & {Solved {\scriptsize(240)}} & {Gap {\scriptsize(43)}} & {Time {\scriptsize(89)}} & {Pace {\scriptsize(124)}} & {Nodes {\scriptsize(192)}} \\ \hline
        {\texttt{ML-Q-RF}} & 197 & 0.017904 & 71.822927 & 0.348300 & 95.867561\\
        \optimal & 198 & 0.013717 & 58.412461 & 0.275232 & 84.902479\\
        \bestruleOPT & 199 & 0.018026 & 61.731941 & 0.284607 & 83.850434\\
        \bestvarQRFtie{0} & 196 & 0.030568 & 70.742677 & 0.362371 & 79.914193\\
        \bestvarOPTtie{0} & 196 & 0.031836 & 72.265518 & 0.365634 & 79.723765\\
    \bottomrule
    \end{tabular}}
\end{table}


Table~\ref{tab:dynamic-sb00} presents the performance of \bestvarQRFtie{0} and \bestvarOPTtie{0}. The results show that, once again, the potential of the dynamic learning setting is inferior to that of the static one. Yet, before drawing any conclusions, it is convenient to check if this behavior may be driven by the zero tolerance considered for ties.

\begin{table}[!htpb]
    \caption{Assessing the potential of \bestvar with different thresholds on ties.}
    \label{tab:dynamic-sbthresholds}
    \centering
    {\addtolength{\tabcolsep}{-2pt} 
    \begin{tabular}{l| r S[round-precision=3] S S S}
    \toprule
        {} & {Solved {\scriptsize(240)}} & {Gap {\scriptsize(46)}} & {Time {\scriptsize(90)}} & {Pace {\scriptsize(129)}} & {Nodes {\scriptsize(188)}} \\ \hline
        \bestvarOPTtie{0} & 196 & 0.025404 & 51.892524 & 0.330583 & 74.003216\\
        \bestvarOPTtie{0.005} & 193 & 0.023866 & 53.646018 & 0.313009 & 74.060050\\
        \bestvarOPTtie{0.01} & 192 & 0.022702 & 53.581657 & 0.314598 & 74.886432\\
        \bestvarOPTtie{0.02} & 192 & 0.023472 & 54.991663 & 0.319552 & 76.420545\\
        \bestvarOPTtie{0.05} & 194 & 0.022263 & 59.864879 & 0.356498 & 77.998123\\
    \bottomrule
    \end{tabular}}
\end{table}


Table~\ref{tab:dynamic-sbthresholds} represents a comparison between five difference values for the tolerance on ties, ranging from $0$ to $0.05$. Although all 5 versions of \bestvar perform comparably well, \bestvarOPTtie{0.01} delivers some improvement in terms of Gap, without compromising much on Time and Pace. Since the performance according to Gap was probably the weakest point of \bestvarOPTtie{0} in Table~\ref{tab:dynamic-sb00}, we present a last table comparing the performance of \bestvarOPTtie{0.01} with respect to \mlrule, \optimal, and \bestrule.

\begin{table}[!htpb]
    \caption{Comparing the potential of \optimal with the best configurations of \bestrule and \bestvar.}
    \label{tab:dynamic-final}
    \centering
    {\addtolength{\tabcolsep}{-2pt} 
    \begin{tabular}{l| r S[round-precision=3] S S S S}
    \toprule
        {} & {Solved {\scriptsize(240)}} & {Gap {\scriptsize(44)}} & {Time {\scriptsize(90)}} & {Pace {\scriptsize(125)}} & {Nodes {\scriptsize(190)}} & {Mean nodes {\scriptsize(190)}} \\ \hline
        {\texttt{ML-Q-RF}} & 197 & 0.016768 & 69.073893 & 0.337245 & 87.813280 & 835.884211\\
        \optimal & 198 & 0.012924 & 56.305667 & 0.266998 & 79.328757 & 672.484211\\
        \bestruleOPT & 199 & 0.016880 & 59.660463 & 0.276658 & 79.091692 & 803.021053\\
        \bestvarOPTtie{0.01} & 192 & 0.026164 & 72.367954 & 0.338298 & 76.221617 & 844.894737 \\
    \bottomrule
    \end{tabular}}
\end{table}


The results in Table~\ref{tab:dynamic-final} confirm that, not only \bestvar is inferior to \optimal, but it is also clearly inferior to \bestrule. We believe that this finding is important on its own. In the context of MILP problems, it has been widely documented that strong branching produces the smallest trees, with no other rule becoming close to it \citep{Achterberg2005,gamrath2018}. In order to assess this in our setting, Table~\ref{tab:dynamic-final} contains one additional column, reporting the standard mean of the average number of nodes, which complements the information about the geometric mean in Nodes, a measure less sensitive to outliers. The geometric means reported in Nodes show that, also in our NLP setting, strong branching imitation tends to produce the smallest trees and, yet, its standard mean is the largest one. A deeper examination of the numerical results behind Table~\ref{tab:dynamic-final} confirms that this is precisely because, although \bestrule performs very well on most instances, there is a non-negligible number of them in which it performs poorly. Thus, despite being the best expert according to Nodes, \bestvar falls clearly behind \optimal and \bestrule according to all other metrics. Interestingly, even \mlrule is very competitive with \bestvar. 

It is natural to wonder what may be the mechanism behind the negative results for \bestvar, \emph{i.e.}, strong branching imitation, in our NLP setting. When compared to \bestvar, \bestrule has the advantage of using information about the magnitudes of the violations of RLT defining identities in Eq.~\eqref{eq:RLTidentity}, and also the weights in Eq.~\eqref{eq:sumweightscriteria} can capture information about the evolution of the branching process. For instance, \texttt{Range} and \texttt{Range Rel} take into account how much the bounds of the variables have changed from their values at the root node. This information allows \bestrule to disregard some variables that might then be selected by strong branching based on their ``myopic'' KPI, which just focuses on LB improvement. The violations of the RLT defining identities and the $w(j,J)$ weights convey additional information that may pay off deeper in the tree. In this sense, although \bestvar may be improving the lower bound faster at the current node compared to \optimal and \bestrule, it may not be reducing so well the size of the infeasibilities. 

\begin{figure}[!htpb]
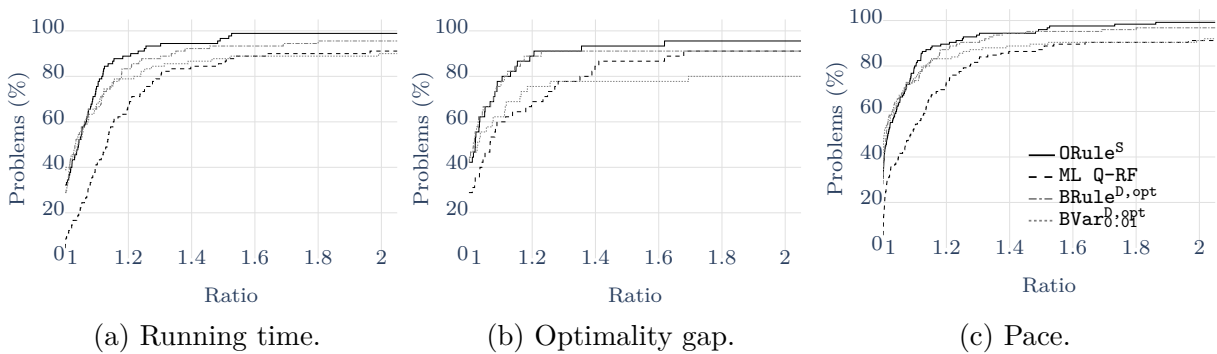

\centering
\begin{subfigure}{0.32\textwidth}
\includesvg[pretex=\scriptsize,width=\textwidth]{images/time2.svg}
    \caption{Running time.}
    \label{subfig:tim2}
\end{subfigure}
\hfill
\begin{subfigure}{0.32\textwidth}
    \includesvg[pretex=\scriptsize,width=\textwidth]{images/gap2.svg}
    \caption{Optimality gap.}
    \label{subfig:gap2}
\end{subfigure}
\hfill
\begin{subfigure}{0.32\textwidth}
    \begin{tikzpicture}
     \node (fig1) at (0,0){\includesvg[pretex=\scriptsize,width=\textwidth]{images/pace2.svg}};
     \node (fig2) at (1,-0.4)
       {\includesvg[pretex=\scriptsize,width=0.4\textwidth]{images/legend2_JGD.svg}};  
    \end{tikzpicture}
    \caption{Pace.}
    \label{subfig:pace2}
\end{subfigure}       
\caption{Performance profiles.}
\label{fig:ppbestvar}
\end{figure}

Figure~\ref{fig:ppbestvar} presents performance profiles comparing \mlrule, \optimal, \bestrule, and \bestvar according to Time, Gap, and Pace. The results are consistent with those reported in Table~\ref{tab:dynamic-final}, but they also allow us to get additional insights. We can see that, for small ratios in the performance profiles, \bestvar is very competitive, and these small ratios already account for a relatively large percentage of the instances. This means that \bestvar is very effective in most of the instances but, as the ratios become larger (\emph{i.e.}, relative differences in performance become larger), it starts falling significantly below \optimal and \bestrule, and, eventually, it falls even below \mlrule. This behavior confirms what we argued above: although \bestrule performs well in most instances, it also performs poorly on a relevant amount of them, which suffices to significantly deteriorate its overall performance.




\section{Conclusions}\label{sec:conclusions}
This paper reveals interesting insights into the effectiveness of different branching rules and variable selection strategies for general non-linear programming problems when using a spatial branch-and-bound algorithm. The analysis was done on RLT for polynomial programming problems however everything can be extended to general non-linear programming. In particular, we analyze three different experts that learning approaches can try to imitate and show the limitations of two of them (\bestvar and \bestrule), which rely on dynamic features compared to the one that is based on static features (\optimal).  These limitations primarily stem from the shortsighted decisions made using local information from the branch-and-bound tree. The dynamic use of features, while seemingly adaptive, often fails to provide long-term benefits, rendering it less effective against the \optimal expert. This underscores the challenge of relying too heavily on information that only pertains to the current state of the tree, without considering broader implications.

Among the three experts, \optimal outperforms the others while \bestrule, which relies on dynamic branching rule selection, performs less effectively because of the myopic nature of the KPI it uses which sets a limit on the attainable outcomes within the dynamic rule framework. Meanwhile, \bestvar is the least effective of the three. Although \bestrule might not always choose the best immediate options, it benefits from incorporating insights into the branching process' progression. This allows \bestrule to exclude certain variables that \bestvar might select based on their immediate, myopic, KPI focused solely on bound improvement. 

The underperformance of \bestvar, which attempts to mimic strong branching, suggests that future efforts in developing branching rules for spatial branching in the context of nonlinear optimization should perhaps not focus solely on approximating strong branching. Instead, these efforts might aim to explore new directions such as extending reinforcement learning to spatial branching which can potentially lead to more effective branching strategies in complex optimization environments. Other directions worth exploring are finding a better KPI to learn in the dynamic setting and enriching the portfolio of branching rules in the static learning setting.

\section*{Acknowledgments}

This work is part of the R\&D projects PID2020-116587GB-I00 and PID2021-124030NB-C32 granted by MICIU/AEI/10.13039/501100011033. This research was also funded by Grupos de Referencia Competitiva ED431C-2021/24 from the Consellería de Cultura, Educación e Universidades, Xunta de Galicia. Ignacio Gómez-Casares acknowledges the support from the Spanish Ministry of Education through FPU grant 20/01555. Bissan Ghaddar's research is supported by the Natural Sciences and Engineering Research Council of Canada Discovery Grant 2017-04185 and by the David G. Burgoyne Faculty Fellowship.

\bibliographystyle{apalike}
\bibliography{references}

\end{document}